\documentclass[12pt]{amsart}
\usepackage{amscd,amssymb}
\usepackage[graph,frame,poly,arc]{xy}  
\usepackage[plainpages,backref,urlcolor=blue]{hyperref}

\theoremstyle{plain}
\newtheorem{thm}[subsection]{Theorem}
\newtheorem{lem}[subsection]{Lemma}

\newtheorem{cor}[subsection]{Corollary}

\theoremstyle{definition}
\newtheorem{rk}[subsection]{Remark}

\newtheorem{ex}[subsection]{Example}

\numberwithin{equation}{section}
\newcommand{\OO}{{\mathcal O}}

\newcommand{\J}{{\mathcal J}}

\newcommand{\A}{{\mathcal A}}

\newcommand{\m}{{\bf m}}
\newcommand{\al}{{\alpha}}
\newcommand{\be}{{\beta}}

\newcommand{\Z}{\mathbb{Z}}
\newcommand{\Q}{\mathbb{Q}}

\newcommand{\C}{\mathbb{C}}
\newcommand{\PP}{\mathbb{P}}

\newcommand{\N}{\mathbb{N}}

\DeclareMathOperator{\Hom}{Hom}
\DeclareMathOperator{\rank}{rank}

\DeclareMathOperator{\coker}{coker}

\begin{document}

\title [On the Milnor monodromy of the complex reflection arrangements]
{On the Milnor monodromy of the irreducible  complex reflection arrangements}

\author[Alexandru Dimca]{Alexandru Dimca$^1$}
\address{Univ. Nice Sophia Antipolis, CNRS,  LJAD, UMR 7351, 06100 Nice, France. }
\email{dimca@unice.fr}

\thanks{$^1$ Partially supported by Institut Universitaire de France.}



\begin{abstract} Using recent results by A. M\u acinic, S. Papadima and  R. Popescu, and a refinement of an older construction of ours, we determine  the monodromy action on $H^1(F(G),\C)$, where $F(G)$ denotes the Milnor fiber of a hyperplane arrangement associated to an irreducible complex reflection group $G$.
\end{abstract}
 
\maketitle


\section{Introduction} \label{sec:intro}

Let $G \subset GL_n(\C)$ be a finite complex reflection group and denote by $\A(G)$ the union of all the reflecting hyperplerplanes of $G$, i.e. the hyperplanes in $\C^n$ fixed by some $g \in G$, $g \ne Id$. Then the complement $M(\A(G))=\C^n \setminus \A(G)$ of the hyperplane arrangement $\A(G)$ has a very special and interesting topology,
namely it is a $K(\pi,1)$-space, see \cite{Bessis} for the long story of the proof of this result.
Since $\A(G)$ is a central hyperplane arrangement, it has a defining equation $f=0$, where $f$ 
is a homogeneous polynomial. One can associate to this setting the Milnor fiber of the arrangement $F(G):f=1$, which is a smooth hypersurface in $\C^n$ endowed with a monodromy morphism $h: F(G) \to F(G)$, see next section for the precise definitions.
The study of the induced monodromy operators 
$$h^j(G): H^j(F(G),\C) \to H^j(F(G),\C)$$
in the case when $G$ is a real reflection group, i.e. a Coxeter or  a Weyl group, was started by S. Settepanella in \cite{Se1}, \cite{Se2}. The easy case of reflection arrangements of rank at most two, and thus corresponding to Milnor fibers of isolated singularities,  is discussed in 
\cite[Section 5]{DL2}.  The case of the reducible complex reflection groups can be reduced to the case $G$ irreducible via Theorem 1.4 (i) in \cite{DNago}.

When $G$ is an irreducible complex reflection group, already to determine the first possibly non-trivial monodromy operator $h^1(G): H^1(F(G),\C) \to H^1(F(G),\C)$ is a challenge. In a recent preprint \cite{MPP}, A. M\u acinic, S. Papadima and R. Popescu have obtained a nearly complete control on the eigenvalues of
the monodromy operator $h^1(G)$ of order $p^s$, with $p$ a prime number and $s$ a positive integer. Some of their main results are stated below, see Theorems \ref{thmMPP}, \ref{thmMPP2} and \ref{thmMPP3}, in order to better understand the contribution of our note.

The irreducible non-exceptional complex reflection arrangements of rank $\geq 3$ consist of the monomial arrangements
$\A(m,m, n)$ with $(m \geq 2,  n \geq 3)$ or $(m=1, n \geq 4)$ (which are in fact the braid arrangements $A_{n-1}$), given as  central hyperplane arrangements in $\C^n$ with coordinates $x_0,...,x_{n-1}$
by
$$\A(m,m,n): f=\prod_{0\leq i<j \leq n-1}(x_i^m-x_j^m)=0$$
and the full monomial arrangements $\A(m,1,n)$ with $m\geq 2, n \geq 3$ given by
$$\A(m,1,n): f=\prod_{k=0,n-1}x_k \prod_{0\leq i<j \leq n-1}(x_i^m-x_j^m)=0.$$

Denote by $F(m,m,n)$ the Milnor fiber of the monomial arrangement $\A(m,m,n)$, and recall
the following results  proved in \cite{MPP}.

\begin{thm}
\label{thmMPP}
For $n=3$, the monodromy operator
$$
h^1: H^1(F(m,m,3),\C) \to H^1(F(m,m,3),\C)
$$
has eigenvalues of order $p^s$ if and only if $p=3$. Moreover, for $p=3$, if we denote the multiplicity of such an eigenvalue by $e^s_3(m)$, then $e^s_3(m) \leq 2$ if $m$ is divisible by $3$, and $e^s_3(m)\leq 1$
otherwise. For $s=1$, both inequalities become equalities.
\end{thm}

\begin{thm}
\label{thmMPP2}
For $n>3$, the monodromy operator
$$
h^1: H^1(F(m,m,n),\C) \to H^1(F(m,m,n),\C)
$$
has eigenvalues of order $p^s$ if and only if $p=3$ and $n=4$. When $n=4$ and $p=3$,  then $e^s_3(m) \leq 1$. For $s=1$, this inequality becomes an  equality.
\end{thm}

Denote by $F(m,1,n)$ the Milnor fiber of the full monomial arrangement $\A(m,1,n)$, and recall
the following result  proved in \cite{MPP}.

\begin{thm}
\label{thmMPP3}
The monodromy operator
$$
h^1: H^1(F(m,1,n),\C) \to H^1(F(m,1,n),\C)
$$
has eigenvalues of order $p^s$ if and only if $p=3$, $n=3$ and $m \equiv 1$ mod $3$. Moreover, for $p=3$ , $n=3$ and $m \equiv 1$ mod $3$,  then $e^s_3(m) \leq 1$. For $s=1$, this inequality is an equality.
\end{thm}

Building on these results by A. M\u acinic, S. Papadima and R. Popescu, and adding a  different approach, in this note we determine completely the eigenvalues of $h^1(G)$ for all irreducible complex reflection groups $G \ne G_{31}$.

\begin{thm}
\label{thmmain}
For $n=3$, the monodromy operator of the monomial arrangement
$$
h^1: H^1(F(m,m,3),\C) \to H^1(F(m,m,3),\C)
$$
has as eigenvalues only cubic roots of unity. Hence $\Delta_V(t)=(t-1)^{3m-1}(t^2+t+1)^2$
if $m$ is divisible by $3$ and $\Delta_V(t)=(t-1)^{3m-1}(t^2+t+1)$ otherwise.

\end{thm}

\begin{thm}
\label{thmmain2}
For $n>3$, the monodromy operator of the monomial arrangement
$$
h^1: H^1(F(m,m,n),\C) \to H^1(F(m,m,n),\C)
$$
has as eigenvalues only cubic roots of unity. Hence 
$$\det (t\cdot Id -h^1|H^1(F(m,m,n),\C))=(t-1)^{qm-1}(t^2+t+1)$$
if $n=4$ and  $h^1=Id$
for $n>4$, where $q= {n \choose 2}.$

\end{thm}

\begin{thm}
\label{thmmain3}
The monodromy operator of the full monomial arrangement 
$$
h^1: H^1(F(m,1,n),\C) \to H^1(F(m,1,n),\C)
$$
has as eigenvalues only cubic roots of unity. Hence 
$$\det (t\cdot Id -h^1|H^1(F(m,1,n),\C))=(t-1)^{qm+n-1}(t^2+t+1)$$
if $n=3$ and  $m\equiv 1$ mod $3$, and $h^1=Id$
otherwise, where $q= {n \choose 2}.$

\end{thm}

The exceptional complex reflection groups of rank $\geq3$ are usually denoted by $G_j$, with $23 \leq j \leq 37$. More precisely, the group $G_j$ has rank 3 for $23 \leq j \leq 27$, rank 4 for $28 \leq j \leq 32$, 
$\rank (G_{33})=5$, $\rank (G_{34})=\rank (G_{35}=E_6)=6$, $\rank (G_{36}=E_7)=7$ and  $\rank (G_{37}=E_8)=8$. For a general reference on complex reflection groups, see \cite{LT}.

\begin{thm}
\label{thmmain4}

The  monodromy operator
$$h^1:H^1(F(G_j),\C) \to H^1(F(G_j),\C)$$
for the exceptional complex reflection group $G_j$ is trivial,
unless  $j=25$. Moreover, $\A(G_{25})$ corresponds to the Hessian arrangement, and one has
$$\det (t\cdot Id -h^1|H^1(F(G_{25}),\C))=(t-1)^{9}(t^4-1)^2.$$
\end{thm}
Here we prove this result for $j \ne 31$, and explain why the case $j=31$ is special. The proof for this latter case uses completely different ideas, and a substantial amount of computer aided computation, so it is presented in a separate paper, see \cite{DSt31}.

\bigskip

The author would like to thank G. Lehrer, A. M\u acinic and S. Papadima for useful discussions.

\section{A brief presentation of our approach}

To describe our new approach, we work in a rather different setting as follows.
Let $V:f=0$ be a hypersurface  of degree $d\geq 3$  in the complex projective space $\PP^n$, defined by a homogeneous polynomial $f \in S=\C[x_0,...,x_n]$. We assume in this paper that $V$ has only isolated singularities and $n \geq 2$.
Consider the corresponding complement $U=\PP^{n}\setminus V$, and the global Milnor fiber $F$ defined by $f(x_0,...,x_n)=1$ in $\C^{n+1}$ with monodromy action $h:F \to F$, $h(x_0,...,x_n)=\exp(2\pi i/d)\cdot (x_0,...,x_n)$. It is known that $H^j(F,\C)=0$ for $1 \leq j <n-1$, see \cite{D1}.

One can consider  the characteristic polynomials of the monodromy, namely
\begin{equation} 
\label{Delta}
\Delta^j_V(t)=\det (t\cdot Id -h^j|H^j(F,\C)),
\end{equation} 
for $j=0,n-1$ and $n$. It is clear that $\Delta^0_V(t)=t-1$, and moreover
\begin{equation} 
\label{Euler}
\Delta^0_V(t)\Delta^{n-1}_V(t)^{(-1)^{n-1}}\Delta^n_V(t)^{(-1)^{n}}=(t^d-1)^{\chi(U)},
\end{equation} 
where $\chi(U)$ denotes the Euler characteristic of the complement $U$, see for instance \cite[Proposition 4.1.21]{D1}.
It follows that the polynomial $\Delta_V(t)=\Delta^{n-1}_V(t)$, also called the Alexander polynomial of $V$, see \cite{L1}, \cite{R}, determines the remaining polynomial $\Delta^n_V(t)$. To find the Alexander polynomial $\Delta_V(t)$, or equivalently, the eigenvalues of the monodromy operator
\begin{equation} 
\label{mono1}
h^{n-1}: H^{n-1}(F,\C) \to H^{n-1}(F,\C)
\end{equation} 
starting from $V$ or  $f$  is a rather difficult problem, going back to O. Zariski and attracting an extensive literature, see for instance \cite{HE}, \cite{L1}, \cite{L2}, \cite{OkaS}, \cite{K},  \cite{D1} for the case $n=2$. 
In this paper we take a new look at a method to determine the Alexander polynomial $\Delta_V(t)$ introduced in \cite{DDuke} and developed in \cite[Chapter 6]{D1}. It is based on the interplay between the Hodge filtration and the pole order filtration on the cohomology groups $H^*(U,\C)$ obtained in \cite{DelDi},
see section 2 below for a brief survey.

The new ingredient put forward in this paper is a careful localization at an isolated weighted homogeneous singularity $(V,p_i)$ of the hypersurface $V$, see section 3. We show that even when the local equation $g'_i=0$ of $(V,p_i)$, obtained say from the global equation $f=0$ of $V$ by choosing an affine chart at $p_i$, is not a weighted homogeneous polynomial, we can control the localizations of rational differential forms as well as when $g'_i$ is replaced by an equivalent  weighted homogeneous polynomial $g_i$  (such a polynomial $g_i$ exists by our hypothesis that $(V,p_i)$ is weighted homogeneous), see Corollaries \ref{corIS} and \ref{corISE}.

Using this localization and its compatibility with the Thom-Sebastiani suspension, we state and prove the main result, which is Theorem \ref{thmvan}. It
gives an effective criterion that a  $d$-th root of unity $\lambda$ is not an eigenvalue for an extension of the monodromy action on the Thom-Sebastiani suspension $V'$ of the hypersurface $V$.
When we start with a plane curve $V$, i.e. for $n=2$, Corollary \ref{corvan2} says that this criterion, applied to $\lambda$ and to its conjugate $\overline {\lambda}$, is enough to conclude that $\lambda$ is not a root of the Alexander polynomial $\Delta_V(t)$.

When a projective hypersurface $V:f=0$ in $\PP^n$ has a singular locus of dimension $\sigma >0$, then we can consider the hypersurface $W=V \cap E: f_W=0$ in $\PP(E)=\PP^{n-\sigma}$,
obtained by intersecting $V$ with a generic codimension $\sigma$ linear subspace $E$. Then $W$ has only isolated singularities, and the cohomology of the Milnor fibers of $f$ and $f_W$ coincide up to the degree $j=n -\sigma-1$, including the monodromy action, see \cite[Theorem 4.1.24]{D1}.
Therefore our results may yield interesting information when $\sigma>0$ as well, as we see below in the case of reflection groups of rank $>3$.

\section{Reformulation of the problem in terms of hypersurface complements}

To fix our convention in a compatible way with \cite{BDS}, \cite{DL1} and \cite{DStCompM}, note that the action of the monodromy on a cohomology class $[\omega] \in H^j(F,\C) $ is given by the formula
\begin{equation} 
\label{mono2}
h^j([\omega])=[(h^{-1})^*(\omega)].
\end{equation} 
Let $\theta=\exp(2\pi i/d)$, and denote by $H^j(F,\C)_{\theta^k}$ the eigenspace of $h^j$ corresponding to the eigenvalue $\theta^k$, for $k=0,1,...,d-1$. Then it is known that we have a natural isomorphism $H^j(F,\C)_1=H^j(U,\C)$ for any $j$. To describe the eigenspaces $H^j(F,\C)_{\theta^k}$ for $0<k<d$, one proceeds as follows. Consider the new homogeneous polynomial
$$f'(x_0,...,x_n,t)=f(x_0,...,x_n)+t^d,$$
called the Thom-Sebastiani $d$-th suspension of $f$, and note that the  hypersurface $V':f'=0$ in $\PP^{n+1}$ is a singular compactification of $F$. More precisely, let $ H$ be the hyperplane in $\PP^{n+1}$ given by $t=0$. Then one has the natural identifications
\begin{equation} 
\label{ident1}
V' \cap H=V \text { and } V' \setminus (V' \cap H)=F.
\end{equation} 
The multiplicative group $\mu_d$ of $d$-th roots of unity acts on $\PP^{n+1}$ via the formula
\begin{equation} 
\label{action1}
\theta \cdot (x_0:...:x_n:t)=(\theta x_0:...:\theta x_n:t)=(x_0:...:x_n:\theta^{-1}t).
\end{equation} 
Then the induced action on $F$ is just the monodromy action, and the action on $V' \cap H=V$ is trivial. The associated exact sequence in the cohomology with compact supports yields
\begin{equation} 
\label{Gysin1}
... \to H^{n}(V,\C) \to H^{n+1}_c(F,\C) \to H^{n+1}(V',\C) \to H^{n+1}(V,\C) \to....
\end{equation} 
and hence, for any nontrivial character $\xi$ of $\mu_d$ we have
\begin{equation} 
\label{char1}
H^{n-1}(F,\C)_{\xi}= H^{n+1}_c(F,\C)_{-\xi}= H^{n+1}(V',\C)_{-\xi}
\end{equation} 
in view of the Poincar\' e duality between $H^{n-1}(F,\C)$ and $H^{n+1}_c(F,\C)$.
Moreover $H^{n+1}(V',\C)_{-\xi}=H^{n+1}_0(V',\C)_{-\xi}$, where 
$$H^m_0(V', \C)=\coker \{ H^m(\PP^{n+1},\C) \to H^m(V',\C)\}$$
denotes the primitive cohomology of $V'$, since the $\mu_d$-action on the cohomology of $ \PP^{n+1}$ is trivial.
See also \cite[Theorem 1.1]{DL1}. 

Note also the convention $H^{n-1}(F,\C)_{\xi}=H^{n-1}(F,\C)_{\theta^k}$ if and only if $\xi=\xi_k:=k$ in the character group $\Z/d\Z=\Hom(\mu_d, \C^*)$.
Set $U'= \PP^{n+1} \setminus V'$, let $p_i\in \PP^n$ for $i=1,...,r$ be the singular points of $V$, and let $p'_i=(p_i:0) \in \PP^{n+1}$ be the corresponding singular points of $V'$. Note that all these points are fixed points under the $\mu_d$-action on $\PP^{n+1}$. For each $i=1,...,r$ choose a Milnor ball $B'_i$ for the isolated hypersurface germ $(V',p'_i)$ (i.e. a small open ball $B'_i$ centered at $p'_i$ such that $V' \cap B'_i$ is topologically a cone with vertex $p'_i$), which in addition is $\mu_d$-invariant. One has, for any integer $s \geq 1$,  a $\mu_d$-equivariant exact sequence
\begin{equation} 
\label{es1}
F^sH^{n+1}(U',\C) \xrightarrow{\rho}  \bigoplus_{i=1,r}F^sH^{n+1}(B'_i \setminus X'_i) \to F^{s-1}H^{n+1}_0(V',\C)\to 0
\end{equation} 
where $F^*$ denotes the Hodge filtration and $X'_i=B'_i \cap V'$, see \cite[(6.3.15)]{D1}. Here the  Hodge filtration on
$H^{n+1}(B'_i \setminus X'_i)$ comes from the identification 
\begin{equation} 
\label{HF1}
F^sH^{n+1}(B'_i \setminus X'_i)= F^{s-1}H^{n+1}_{\{p'_i\}}(V'), 
\end{equation} 
see \cite{DurSai},  \cite{DurHain}, \cite{Nav} or pp. 200-201 in \cite{D1}. Moreover, the morphism $\rho$ is induced by the restriction of rational differential forms, and this gives an additional exact sequence 
\begin{equation} 
\label{es2}
S_{(n-s+2)d-n-2} \xrightarrow{\rho'}  \bigoplus_{i=1,r}F^sH^{n+1}(B'_i \setminus X'_i) \to F^{s-1}H^{n+1}_0(V',\C)\to 0
\end{equation} 
when all the singularities $(V,p_i)$ are weighted homogeneous, see \cite[(6.3.16)]{D1}, but note that filtration $F^{s-1}$ on the first cohomology group there should be replaced by $F^s$. The morphism $\rho'$ is given by the following formula $\rho'(h)=\rho(\omega(h))$, where 
\begin{equation} 
\label{forms1}
\omega(h)=\frac{h\Omega}{(f')^{n+2-s}}
\end{equation} 
and $\Omega$ is the contraction by the Euler vector field of the top differential form 
$$dx_0 \wedge ...\wedge dx_n \wedge dt.$$
\begin{rk} 
\label{rknoniso}
This method to study the topology of projective hypersurfaces was extended to hypersurfaces with non-isolated singularities in \cite{HK}, and this may open the door to new applications, e.g. the computation of the monodromy $h^2$ in the case of a plane arrangement in $\PP^3$.
\end{rk} 

\section{On the local cohomology of an isolated weighted homogeneous singularity}

\subsection{The absolute case}

Let $\OO_0$ be the ring of complex analytic function germs at the origin of $\C^n$ and denote by $\m$ its unique maximal ideal. Let $g \in \m$ and assume that $X_g:g=0$ is an isolated singularity.
Fix $B$ a Milnor ball for $X_g$ such that, in particular, the germ $g$ has a representative, also denoted by $g$, defined on $B$.  The set of such balls $B$ forms clearly a projective system, and hence we can define 
\begin{equation} 
\label{coho1}
H^n(B_g\setminus X_g,\C)=\varinjlim H^n(B\setminus X_g,\C),
\end{equation} 
an injective limit where all the morphisms are isomorphisms.
Note that in fact one has a natural isomorphism $H^n(B_{g}\setminus X_{g},\C)=H^n_{\{0\}}(X_{g},\C)$, the local cohomology of $X_g$ with support at the origin, so there is already an intrinsic notation for this object. However, in view of the shift in the Hodge filtration \eqref{HF1} and the next construction, we prefer our ad-hoc definition. 

Let $\Omega^n_0$ denote the $\OO_0$-module of germs of holomorphic $n$-forms at the origin of $\C^n$. Then there is a well defined $\C$-linear map
\begin{equation} 
\label{coho2}
L_g:  \Omega^n_0 \to   H^n(B_g\setminus X_g,\C),
\end{equation} 
sending $\omega \in \Omega^n_0$ to the cohomology class of the meromorphic form $\frac{\omega}{g}$ in some cohomology group $H^n(B\setminus X_g,\C)$ such that the germs $\omega $ and $g$  are defined on the Milnor ball $B$ (and then regard this cohomology class in the limit $H^n(B_g\setminus X_g,\C)$ in the obvious way). This construction is natural in the obvious sense: if $\phi:(\C^n,0) \to (\C^n,0)$ is an analytic isomorphism germ, and if we set $g'=\phi^*(g)=g \circ \phi$, then 
\begin{equation} 
\label{funct1}
\phi^*(L_g(\omega))=L_{g'}(\phi^*(\omega)).
\end{equation} 
Assume now that $g$ is a weighted homogeneous polynomial of type $(w_1,...,w_n;e)$, with $w_j$ strictly positive integers, having an isolated singularity at the origin. Denote by $y_1,...,y_n$ the coordinates on $\C^n$ and let $(y^{\al})_{\al \in A}$ be the monomial basis of the space of weighted homogeneous polynomials in $\C[y_1,...,y_n]$ of degree $e-w_1-...-w_n$. In this case we have a canonical isomorphism 
\begin{equation} 
\label{coho3}
 H^n(B_g\setminus X_g,\C)=H^n(U_g,\C),
\end{equation} 
where $U_g=\C^n\setminus g^{-1}(0)$, induced by the inclusions $B\setminus X_g \to U_g$. Using this isomorphism, we define $F^nH^n(B_g\setminus X_g,\C)$ to be the subspace in $H^n(B_g\setminus X_g,\C)$ corresponding to the Hodge filtration subspace $F^nH^n(U_g,\C)$ in $H^n(U_g,\C)$. It is known that the cohomology group $H^n(U_g,\Q)$ has a pure Hodge structure of weight $n+2$, see \cite{D1}, p. 203.
If we set $\omega_n=dy_1 \wedge ....\wedge dy_n$, then it is known that the cohomology classes
\begin{equation} 
\label{coho4}
\epsilon_{\al}=\frac{y^{\al}\omega_n}{g}
\end{equation}
for $\al \in A$ give a basis for  $F^nH^n(U_g,\C)$, see pp. 202-203 in \cite{D1}.
We denote by $\epsilon_{\al}'$ the element in  $F^nH^n(B_g\setminus X_g,\C)$  corresponding to $\epsilon_{\al}$. Then we have the following result.

\begin{thm}
\label{thmIS}
With the above notation, the following hold.

\begin{enumerate}

\item If $\omega \in \Omega^n_0$ is given by $\omega= \sum_{\be \in \N^n}c_{\be}y^{\be}\omega_n$, then 
$$L_g(\omega)=\sum_{\al \in A}c_{\al}\epsilon'_{\al}.$$
In particular, the image of $L_g$ is exactly $F^nH^m(B_g\setminus X_g,\C)$
\item Let $a(g)=\min \{s \in \N \ : \  sw_j>e-w_1-...-w_n \text{ for all } j \}$. Then one has
$\m^{a(g)} \Omega^n_0 \subset \ker L_g$, and hence there is an induced surjective linear map
$$L_g:\Omega^n_0/\m^{a(g)} \Omega^n_0 \to F^nH^n(B_g\setminus X_g,\C)=F^{n-1}H^m_{\{0\}}(X_{g},\C).$$

\end{enumerate}

\end{thm}

\proof It is enough to show that if $h(y)$ is a weighted homogeneous polynomial of degree
$e'\ne e_0=e-w_1-...-w_n$, then the cohomology class of $\frac{h(y)\omega_n}{g}$ in $H^n(U_g,\C)$ is trivial. Consider the morphism $h_t: U_g \to U_g$, given by 
$$h_t(y)=(t^{w_1}y_1,...,t^{w_n}y_n),$$
for $t \in \C^*$. Since $h_t$ is clearly homotopic to the identity $h_1$, it follows that the induced morphism in cohomology is the identity. However we get
$$h_t^*([ \frac{h(y)\omega_m}{g}])=t^{e'-e_0}[\frac{h(y)\omega_m}{g}],$$
for any $t \in \C^*$, and this implies that the cohomology class $[\frac{h(y)\omega_m}{g}]$ should vanish.

\endproof

For the next claim, use the functoriality described in \eqref{funct1}, and the obvious fact that
$\phi^*$ preserves $\m^{a(g)} \Omega^n_0$ and the Hodge filtration on the local cohomology groups $H^n_{\{0\}}(X_{g},\C)$, see \cite{DurSai},  \cite{DurHain}, \cite{Nav}.

\begin{cor} 
\label{corIS} 
Let $g$ be a weighted homogeneous polynomial of type $(w_1,...,w_n;e)$, having an isolated singularity at the origin. Let $\phi:(\C^n,0) \to (\C^n,0)$ be an analytic isomorphism germ, and set $g'=\phi^*(g)=g \circ \phi$. Then there is an induced surjective linear map
$$L_{g'}:\Omega^n_0/\m^{a(g)} \Omega^n_0 \to F^nH^n(B_{g'}\setminus X_{g'},\C)=F^{n-1}H^n_{\{0\}}(X_{g'},\C),$$
where the integer $a(g)$ is defined in Theorem \ref{thmIS}(2).

\end{cor}

\subsection{The relative case and Thom-Sebastiani construction}

Let $g \in \m$ and assume that $X_g:g=0$ is an isolated singularity. Consider the suspension 
$$G(y,t)=g(y)+t^d$$
as a germ at the origin of $\C^{n+1}$, with coordinates $y_1,...,y_n,t$. Then the corresponding cohomology group $H^{n+1}(B_G \setminus X_G,\C)$ has a natural $\mu_d$-action, coming from the $\mu_d$-action on $\C^{n+1}$ defined by
$$\lambda \cdot (y,t)=(y, \lambda ^{-1}t),$$
for any $\lambda \in \mu_d$. Note that a germ of a holomorphic $(n+1)$-form $\omega$ at the origin of $\C^{n+1}$ can be written as 
$$ \omega= \sum_{m \in \N}
 \omega ^{(m)}\wedge t^mdt,$$
 where $\omega ^{(m)} \in \Omega^n_0$. Using this decomposition, the linear map $L_G$ defined in \eqref{coho2} can be refined taking into account the isotypical components with respect to the $\mu_d$-action, and we get, for $k=1,2,...,d$
\begin{equation} 
\label{coho2E}
L_G:  \Omega^n_0 \wedge t^{k-1}\C\{t^d\}dt  \to   H^{n+1}(B_G\setminus X_G,\C)_{\xi_k},
\end{equation} 
where $\xi_d=\xi_0$ is the trivial character. This construction is clearly functorial as in \eqref{funct1} with respect to germs of isomorphisms $\Phi: (\C^{n+1},0) \to (\C^{n+1},0)$,
where $\Phi$ is a product $\phi \times Id_{\C}$, with $\phi: (\C^{n},0) \to (\C^{n},0)$ an isomorphism as in \eqref{funct1}.

Assume now that $g$ is a weighted homogeneous polynomial of type $(w_1,...,w_n;e)$, with $w_j$ strictly positive integers, having an isolated singularity at the origin. Let $\gamma(e,d)$ be the greatest common divisor of $e,d$ and $\mu(e,d)= ed/ \gamma(e,d)$ be their least common multiple.

Then the suspention $G$ is weighted homogeneous of type
$(d_1w_1,...,d_1w_n, e_1; \mu(e,d))$, where $d_1=d/\gamma(e,d)$ and $e_1=e/\gamma(e,d)$. Let $(y^{\al})_{\al \in A_k}$ be the monomial basis of the space of weighted homogeneous polynomials in $\C[y_1,...,y_n]$ of degree $\mu(e,d)-d_1w_1-...-d_1w_n-e_1k$ with respect to the new weights $wt(y_j)=d_1w_j$, for $k=1,2,...,d$. Then it is clear that the cohomology classes
\begin{equation} 
\label{coho4E}
\epsilon_{\al}=\frac{y^{\al}t^{k-1}\omega_n\wedge dt }{G}
\end{equation}
for $\al \in A_k$ give a basis for  $F^{n+1}H^{n+1}(U_G,\C)_{\xi_k}.$
The corresponding basis in the vector space $F^{n+1}H^{n+1}(B_G\setminus X_G,\C)_{\xi_k}$ is denoted by 
$\epsilon'_{\al}$. We have the following result, with exactly the same proof as for Theorem \ref{thmIS}.

\begin{thm}
\label{thmISE}
With the above notation, the following hold.

\begin{enumerate}

\item If $\omega \in \Omega^n_0 \wedge t^{k-1}\C\{t^d\}dt$ is given by $\omega= \sum_{\be \in \N^n, m \in \N}c_{\be,m}y^{\be} t^{k-1+m} \omega_n \wedge dt$, then 
$$L_G(\omega)=\sum_{\al \in A_k}c_{\al}\epsilon'_{\al}.$$
In particular, the image of $L_G$ from \eqref{coho2E} is exactly $F^{n+1}H^{n+1}(B_G\setminus X_G,\C)_{\xi_k}$.
\item Let $a(g,k)=\min \{s \in \N \ : \ sd_1w_j>\mu(e,d)-d_1w_1-...-d_1w_n-e_1k \text{ for all } j \}$. Then 
$\m^{a(g,k)} \Omega^n_0 \wedge t^{k-1}\C\{t^d\}dt \subset \ker L_G$, and hence there is an induced surjective linear map
$$L_G:(\Omega^n_0/\m^{a(g,k)} \Omega^n_0) \wedge t^{k-1}dt \to F^{n+1}H^{n+1}(B_G\setminus X_G,\C)_{\xi_k}=F^{n}H^{n+1}_{\{0\}}(X_{G},\C)_{\xi_k}.$$

\end{enumerate}

\end{thm}

We also have the following version of Corollary \ref{corIS}.
\begin{cor} 
\label{corISE} 
Let $g$ be a weighted homogeneous polynomial of type $(w_1,...,w_n;e)$, having an isolated singularity at the origin. Let $\phi:(\C^n,0) \to (\C^n,0)$ be an analytic isomorphism germ, and set $g'=\phi^*(g)=g \circ \phi$. Let $G'(y,t)=g'(y)+t^d$.
Then there is an induced surjective linear map
$$L_{G'}:(\Omega^n_0/\m^{a(g,k)} \Omega^n_0) \wedge t^{k-1}dt \to F^{n+1}H^{n+1}(B_{G'}\setminus X_{G‘},\C)_{\xi_k}=F^{n}H^{n+1}_{\{0\}}(X_{G'},\C)_{\xi_k}.$$
where the integer $a(g,k)$ is defined in Theorem \ref{thmISE}(2).

\end{cor}

\begin{ex} 
\label{exsurf1}
(i) Consider the case $n=2$ and $g(y_1,y_2)=y_1^3+y_2^3$. Then $g$ is weighted homogeneous of type $(1,1;3)$. Consider the suspension $G(y_1,y_2,t)=g(y_1,y_2)+t^{3m}$ for some integer $m\geq 1$, which is weighted homogeneous of type $(m,m,1;3m)$. Using the absolute case discussed above, we see that a basis for $F^3H^3(B_G \setminus X_G,\C)$ is given by the following forms. Fix an integer $k \in \{1,2,...,3m\}$ and then fix a monomial basis $(y^{\al})_{\al \in A_k}$ of the space of weighted homogeneous polynomials in $\C[y_1,y_2]$ of degree $d-d_1w_1-...-d_1w_n-k=3m-m-m-k=m-k$, where $wt(y_i)=d_1w_i=m$.  It follows that only the value $k=m$ gives a nontrivial vector space, which 1-dimensional with the basis $1$. It follows that 
$$F^3H^3(B_G \setminus X_G,\C)=F^3H^3(B_G \setminus X_G,\C)_{\xi_m},$$
and $a(g,m)=1.$

(ii) Consider now the case $n=2$ and $g(y_1,y_2)=y_1^m+y_2^m$, for some integer $m \geq 3$. Then $g$ is weighted homogeneous of type $(1,1;m)$. Consider the suspension $G(y_1,y_2,t)=g(y_1,y_2)+t^{qm}$, for some integer $q>0$. Then $G$ is weighted homogeneous of type $(q,q,1;qm)$. Using the absolute case discussed above, we see that a basis for $F^3H^3(B_G \setminus X_G,\C)$ is given by the following forms. Fix an integer $k \in \{1,2,...,qm\}$ and then fix a basis $(y^{\al})_{\al \in A_k}$  of the space of weighted homogeneous polynomials in $\C[y_1,y_2]$ of degree $d-d_1w_1-...-d_1w_n-k=qm-2q-k=q(m-2)-k$, where $wt(y_i)=d_1w_i=q$.  It follows as above that 
only the values $k=qk_1$, for $1 \leq  k_1 \leq m-2$ give a nontrivial vector space and 
that one has
$$F^3H^3(B_G \setminus X_G,\C)=\bigoplus_{1 \leq  k_1 \leq m-2}F^3H^3(B_G \setminus X_G,\C)_{\xi_{3k_1}},$$
 and $a(g,qk_1)=m-k_1-1.$

(iii) Consider again  the case $n=2$ and $g(y_1,y_2)=y_1^m+y_2^m$, for some integer $m \geq 3$.  Consider the suspension $G(y_1,y_2,t)=g(y_1,y_2)+t^{d}$, which is weighted homogeneous of type $(d_1,d_1,m_1;\mu(m,d))$, where $d_1=d/\gamma(m,d)$ and  $m_1=m/\gamma(m,d)$. Using the absolute case discussed above, we see that a basis for $F^3H^3(B_G \setminus X_G,\C)$ is given by the following forms. Fix an integer $k \in \{1,2,...,d\}$ and then fix a basis $(y^{\al})_{\al \in A_k}$  of the space of weighted homogeneous polynomials in $\C[y_1,y_2]$ of degree $\mu(m,d)-2d_1-km_1=d_1(m-2)-km_1$.  Since $m_1$ and $d_1$ are relatively prime, it  follows that only the values $k=k_1d_1$, for $1 \leq  k_1 < \gamma(m,d)$ give a nontrivial vector space.
Such a value for $k$ corresponds to the eigenvalue $ \theta^k=\exp(2\pi i k/d)= \exp(2 \pi i k_1/\gamma(m,d))$. It follows that 
and $a(g,k) =m-1-k_1m_1.$

Note that one may write 
$$a(g,k) =m-1-k_1m_1= \frac{mk'}{d}-1,$$
for $k'=d-k$, a formula needed in our Remark \ref{rkBDS2} below.

\end{ex}

\section{How to prove the vanishing of some monodromy eigenspaces}

\subsection{The general approach}

Consider a projective hypersurface $V$ in $\PP^n$ having only isolated, weighted homogeneous singularities, and recall the exact sequence \eqref{es2} for $s=n+1$. Then, in terms of $\mu_d$-isotypic components we have the following exact sequence
\begin{equation} 
\label{es2E}
S_{d-n-1-k} \xrightarrow{\rho_k}  \bigoplus_{i=1,r}F^{n+1}H^{n+1}(B'_i \setminus X'_i,\C) _{\xi_k}\to F^{n}H^{n+1}_0(V',\C)_{\xi_k}\to 0
\end{equation} 
for any $k=1,...,d-1$. Here we can clearly replace $F^{n+1}H^{n+1}(B'_i \setminus X'_i) _{\xi_k}$
by the more intrinsic object $F^{n+1}H^{n+1}(B_{G'_i} \setminus X_{G'_i}, \C) _{\xi_k}$, where $g'_i=0$ is a (local) equation for the isolated singularity $(V,p_i)$ and $G'_i$ is the $d$-th suspension of $g'_i$ as above. Moreover, the morphism $ \rho_k$ sends a homogeneous polynomial $h(x)$ in $(x_0,...,x_n)$ of degree $d-n-1-k$ to the set of $r$ cohomology classes defined by the restrictions of the rational differential form
\begin{equation} 
\label{forms2}
\frac{h(x)t^{k-1}\Omega}{f'}
\end{equation} 
to the complements $B_{G'_i} \setminus X_{G'_i}$ for $i=1,...,r$. If we choose the hyperplane $H_0: x_0=0$ such that all the singularities $p_i$ are in $\C^n=\PP^n \setminus H_0$, and chose $y_i=x_i/x_0$ as coordinates on $\C^n$, then the above restrictions have the form 
$$\eta_i=\frac{h_i(y)t^{k-1}\omega_n\wedge dt}{g'_i},$$
 for $h_i$ analytic germs at the points $b_i$ in $\C^n$ corresponding to the singularities $p_i$'s. We denote by $\OO_{b_i}$ the ring of such analytic function germs at $b_i$, and by $\m _{b_i} \subset \OO_{b_i}$ the corresponding maximal ideals.

To show that $F^{n}H^{n+1}_0(V',\C)_{\xi_k}=0$ for some $k$, we have to show that the morphism $\rho_k$ is surjective. In view of Corollary \ref{corISE}, this morphism can be factor as follows
\begin{equation} 
\label{factor1}
S_{d-n-1-k} \xrightarrow{u} \bigoplus_{i=1,r}  \OO_{b_i}/ \m _{b_i}^{a(g_i,k)}  \xrightarrow{v} \bigoplus_{i=1,r}F^{n+1}H^{n+1}(B_{G'_i} \setminus X_{G'_i}, \C) _{\xi_k},
\end{equation} 
where $g_i$ is a weighted homogeneous polynomial right equivalent to $g'_i$. Notice that we do not need the actual polynomial $g_i$, just its homogeneity type in order to compute the invariant
$a(g_i,k)$. The first morphism $u$ is an evaluation map. It takes a homogeneous polynomial $h$ to the classes of the germs $h_i(y)$ at the points $b_i$ of the polynomial $h(1,y_1,...,y_n)$, while the morphism $v$ sends the $r$-tuple of germ classes $(h_1(y),...,h_r(y))$ to the $r$-tuple of cohomology classes
$([\eta_1],....,[\eta_r])$ defined above. Moreover, the morphism $v$ is surjective by Corollary \ref{corISE}.

On the other hand, recall the following result, see for instance  Corollary 2.1 in \cite{BeSo}.

\begin{lem}
\label{lemvan}
The evaluation morphism
$u=eval_N:S_N \to  \bigoplus_{i=1,r}  \OO_{b_i}/ \m _{b_i}^{a(g_i,k)}   $
is surjective for any $N \geq \sum_{i=1,r} a(g_i,k)-1$.
\end{lem}
More precisely, for a given $k$, let $I_k \subset [1,d-1]$ be the subset consisting of all the indices $i$ such that $F^{n+1}H^{n+1}(B_{G'_i} \setminus X_{G'_i}, \C) _{\xi_k} \ne 0$. Then the above discussion implies the following result.

\begin{thm}
\label{thmvan}
If  $N=d-n-1-k \geq \sum_{i \in I_k} a(g_i,k)-1$, then $F^{n}H^{n+1}_0(V',\C)_{\xi_k}=0$.
\end{thm}

\subsection{The case $n=2$}

We assume in this subsection that $V$ is a plane curve having only weighted homogeneous singularities. Then each of the singularities $(V',p_i')$ is a weighted homogeneous singularity and hence the corresponding complement $H^3(B_{G_i'} \setminus X_{G'_i}, \C)$ has a pure Hodge structure of weight 5, see for instance \cite{D1} p. 245 and pp. 66-67. It follows that
one has
\begin{equation} 
\label{Hodge1}
H^3(B_{G_i'} \setminus X_{G'_i}, \C)=H^{3,2}(B_{G_i'} \setminus X_{G'_i}, \C) \oplus H^{2,3}(B_{G_i'} \setminus X_{G'_i}, \C),
\end{equation} 
where $H^{3,2}(B_{G_i'} \setminus X_{G'_i}, \C)=F^3H^3(B_{G_i'} \setminus X_{G'_i}, \C)$
 and $H^{2,3}(B_{G_i'} \setminus X_{G'_i}, \C)$ is its complex conjugate.
It follows, via the exact sequence \eqref{es1}, that $H^3(V',\Q)$ is a pure Hodge structure of weight 3 such that 
\begin{equation} 
\label{Hodge2}
H^3(V',\C)=H^{2,1}(V',\C) \oplus H^{1,2}(V',\C),
\end{equation} 
where $H^{2,1}(V',\C)=F^2H^3(V',\C)$ and $H^{1,2}(V',\C)$ is its complex conjugate. Next, for any nontrivial character $\xi_k$ of $\mu_d$, we have similar relations
\begin{equation} 
\label{Hodge3}
H^3(V',\C)_{\xi_k}=H^{2,1}(V',\C)_{\xi_k} \oplus H^{1,2}(V',\C)_{\xi_k},
\end{equation} 
and $\dim H^{2,1}(V',\C)_{\xi_k}= \dim H^{1,2}(V',\C)_{\xi_{d-k}}$. These relations combined with the formula \eqref{char1} yield the following.

\begin{cor} 
\label{corvan2}
Let $V$ be a plane curve of degree $d$ having only weighted homogeneous singularities. Then for any nontrivial character $\xi_k$ of $\mu_d$ one has $H^1(F,\C)_{\xi_{k}}=0$ if and only if $F^2H^3(V',\C)_{\xi_{k}}=F^2H^3(V',\C)_{\xi_{d-k}}=0.$
\end{cor} 

This corollary says that for such plane curves the vanishing of an eigenspace of the Milnor monodromy on $H^1(F,\C)$ can be tested using only rational forms with poles of order 1.
In other cases, poles of higher orders are necessary, see for instance the case of nodal hypersurfaces of dimension $\geq 2$ in \cite[Theorem 6.4.5]{D1}. Our new approach outline here will be extended to such more general situations in a subsequent paper.

\begin{rk} 
\label{rkBDS2}
Let $V$ be a line arrangement in $\PP^2$. Then our Theorem \ref{thmvan} looks very similar to \cite[Theorem 2]{BDS}. In fact, by the formula at the end of Example \ref{exsurf1}, (iii), we see that the ideal
$\J_y^{(>k')}$ from \cite[Theorem 2]{BDS} coincides with the ideal $\m _{b_i}^{a(g_i,k)} $, where $y \in V$ is the singular point corresponding to $b_i$. However, note that the target space for the two evaluation morphisms are distinct, because the corresponding sums involve different sets of singular points. It is not clear at this stage whether one of these two results implies the other.

\end{rk} 

\section{The irreducible complex reflection arrangements of rank $\geq 3$}

For a line arrangement $\A$ in $\PP^2$, it  is a major open question whether the monodromy operator $h^1:H^1(F,\C) \to H^1(F,\C)$ is combinatorially determined, i.e. determined by the intersection lattice $L(\A)$, see \cite{OT}. 
Several interesting examples have been computed by D. Cohen, A. Suciu, A. M\u acinic, S. Papadima, M. Yoshinaga, see \cite{CS}, \cite{S1}, \cite{MP}, \cite{S2}, \cite{Yo0}.
When the line arrangement $\A$ has only double and triple points, then a complete positive answer is given by S. Papadima and A. Suciu in \cite{PS}. However, the determination of the eigenvalues of $h^1$ in general remains a very difficult question. Following \cite{MPP}, we discus below the case of (generic $3$-dimensional sections) of the  irreducible non-exceptional reflection arrangements of rank $\geq3$. 

\medskip

 First we give  the proof of Theorem \ref{thmmain}.
The curve $V$, obtained as the union of all the lines in  $\A(m,m,3)$, has two types of singularities, the triple points, with local equation $y_1^3+y_2^3=0$ and 
 the three $m$-fold intersection points,  $p_1=[1:0:0]$, $p_2=[0:1:0]$ and $p_3=[0:0:1]$, with local equation $y_1^m+y_2^m=0$.

 In this case the degree of the curve $V$ is $d=3m$. Using Theorem \ref{thmMPP}, it remains to show that 
$H^1(F,\C)_{\xi_{k}}=0$ for any integer $k$ with $0 <k<3m$ and $k \ne m$, $k \ne 2m$, where $F=F(m,m,3)$ is the corresponding Milnor fiber.
By Corollary \ref{corvan2}, it is enough to show for such $k$'s, the vanishing $F^2H^3(V',\C)_{\xi_{k}}=0$. To do this we use Theorem \ref{thmvan}. First note that, in view of Example
\ref{exsurf1} (i), the triple points give no contribution (i.e. their indices are not in the set $I_k$).
There are at most 3 singular points left, namely $p_1,p_2,p_3$ when $m \ne 3$. Note that the character $\xi_k$ correspond to the eigenvalue $\theta^k=\exp(2\pi i k/3m)$. If such an eigenvalue occurs, it must be a root of at least one of the local Alexander polynomials of the singularities of $V$, see \cite[Corollary 6.3.29]{D1} or \cite{L1}. The local Alexander polynomial of a triple point has only cubic roots of unity, while the local Alexander polynomials of $(V,p_i)$ for $i=1,2,3$ have only $m$-th roots of unity. It follows that $\theta^k$, supposed not a cubic root of unity, can be an eigenvalue only if $k$ is divisible by $3$, say $k=3k_1$.
Then, by 
Example
\ref{exsurf1} (ii), we know that one has $a(g_i,k)=m-k_1-1$ for $i=1,2,3$. It follows that
$$N=3m-3-k \geq 3(m-k_1-1)-1,$$
and hence our claim is proved by Theorem \ref{thmvan}.

\begin{rk} 
\label{rkMPP}
(i) Note that our approach is complementary to that used in \cite{MPP}. Indeed, it seems not easy to obtain the multiplicities of the cubic roots of unity for the monomial line arrangement by the method described above. 

(ii) Theorem \ref{thmmain} was checked for $2 \leq m \leq 25$ using a completely different, computational point of view in \cite{DStCompM}.
\end{rk} 

\medskip

 We give now the proof of Theorem \ref{thmmain2}. The curve $V$ is in this case obtained by taking first the trace of the monomial arrangement $\A(m,m,n)$ on a $3$-dimensional generic linear subspace $E$ in $\C^n$ with $0 \in E$, and then working in the corresponding projective plane $\PP^2=\PP(E)$. It is easy to see that this curve $V$ has degree $d=qm$. Assume in the sequel that $m \geq 4$ (otherwise the claim of Theorem \ref{thmmain2} is obvious, using the local Alexander polynomials).
 
Then $V$ has  $n_m=q$ ordinary $m$-points, $n_3={n \choose 3} m^2$ triple points and 
$$n_2=\frac{1}{2}{n \choose 2}{n-2 \choose 2} m^2$$
double points. It follows exactly as above, by looking at the local Alexander polynomials, that we have only to discard eigenvalues $\theta^k=\exp(2\pi i k/qm)$ of order $>3$ and a divisor of $m$. It follows that one should have $k=qk_1$ for some integer $k_1$. Then, exactly as above we get $a(g_i,k)=m-k_1-1$ for all the $q$ singular points $g_i=0$ of multiplicity $m$ by Example \ref{exsurf1} (ii). Then
$$N=qm-3-k \geq q(m-k_1-1)-1,$$
and hence our claim is proved by Theorem \ref{thmvan}, since as above the double and the triple points can be ignored in this computation.

Next we consider the case of the full monomial arrangement $\A(m,1,n)$. For $n>3$, proceed as above and take first a generic $3$-dimensional section. Then the corresponding new curve $V_f$ is obtained from the previous curve $V$, coming from a monomial arrangement, by adding $n$ new lines $L_i$ corresponding to $x_i=0$ for $i=0,...,n-1$. Each line, say  $L_0$ to fix the ideas, passes through $(n-1)$ points of multiplicity $m$ of the curve $V$, e.g. $L_0$ passes through the points corresponding to $x_0=x_1=0$, $x_0=x_2=0$,..., $x_0=x_{n-1}=0$. In conclusion, the $q$ points of multiplicity $m$ on $V$ become now $q$ points of multiplicity $m'=m+2$ on $V_f$, the triple points in $V$ stay triple points on $V_f$ and there are some additional nodes on $V_f$. Note that the degree of $V_f$ is $d=qm+n$.

To prove Theorem \ref{thmmain3}, we have to show that there are no eigenvalues $\theta^k=\exp(2\pi i k/d)$ of order $>3$ and a divisor of $m'$. Let $\gamma(m',d)$ be the greatest common divisor of $m'$ and $d$ and set $m'=\gamma(m',d)m'_1$ and $d=\gamma(m',d)d_1$.
If the order of $\theta^k$ is a divisor of $m'$, it follows that $km'=k\gamma(m',d)m'_1$ is divisible by $d=\gamma(m',d)d_1$, and hence $k$ is divisible by $d_1$, say $k=d_1k_1$. 

Since $1 \leq k <d$, it follows that $1 \leq k_1 <\gamma(m',d)$.
Using Example \ref{exsurf1} (iii),
we infer that $a(g,k) = m'-1-k_1m'_1$. Moreover, one has
$$d_1=\frac{d}{\gamma(m',d)}=\frac{qm'-2q+n}{\gamma(m',d)}=qm_1'-\frac{2}{\gamma(m',d)}q+n \geq q(m_1'-1),$$
since ${\gamma(m',d)}>3$, being a multiple of the order of $\theta^k$.
 It follows that
$$(d-n-1-k) - (\sum_{i \in I_k} a(g_i,k)-1) =d_1( {\gamma(m',d)} -k_1)-3-q(m'-1-k_1m'_1)+1 \geq $$
$$\geq q(m'-1)(\gamma(m',d)-k_1-1)+qk_1m'_1-2 > 0,$$
since in the last sum the first term is clearly $\geq 0$, while the second satisfies
$$qk_1m'_1-2 \geq n-2 \geq 1.$$
Again the double and the triple points can be ignored in this computation, and hence the proof of 
Theorem \ref{thmmain3} is complete.

Finally, we  prove Theorem \ref{thmmain4} for $j \ne 31$.
In fact, the results in \cite{MPP}, \cite{OT} completely determine the monodromy operator 
$$h^1:H^1(F(G_j),\C) \to H^1(F(G_j),\C)$$
 for  $j \ne 31$. More precisely, by Theorem 1.2 in  \cite{MPP}, one has $h^1=Id$ for $j \ne 25, 31$ and for $j=25$, which corresponds to the Hessian arrangement, one has
$$\det (t\cdot Id -h^1|H^1(F(G_{25}),\C))=(t-1)^{9}(t^4-1)^2,$$
see for instance \cite[Remark 3.3 (iii)]{BDS}.
Indeed, consider as an example the case of the exceptional group $H_3=G_{23}$. Then Table C.4 in \cite{OT} shows that the associated curve $V(G_{23})$ in $\PP^2$ has degree $d=15$, 15 nodes, 10 triple points and 6 points of multiplicity 5, corresponding to the isotropy group $I_2(5)$.  Looking at the local Alexander polynomials shows that the monodromy eigenvalues should have an order which divides 2,3 or 5.  But eigenvalues of these orders are excluded by Theorem 1.2 in  \cite{MPP}.

\begin{rk} 
\label{rkexce}

In the case of the remaining exceptional group $G_{31}$, Table C.12 in \cite{OT} shows that the corresponding curve $V(G_{31})$ in $\PP^2$ has degree $d=60$, 360 double points, 320 triple points and 30 points of multiplicity 6, corresponding to the isotropy group $G(4,2,2)$. If we try to apply our method to this curve, in order to exclude the roots of unity of order 6, we have to take in the above notation $\gamma(6,d)=6$, $d_1=10$, $k=10$ or $k=50$.
When $k=10$,  for any point $p_i$ of multiplicity 6 we get $a(g_i,10)=4$.
Therefore the inequality in Theorem \ref{thmvan} becomes in this case
$60-3-10\geq 30 \cdot 4-1$, which is clearly false. Hence our method does not work in this case. For 
the computation of the Milnor monodromy of arrangement $\A(G_{31})$ using a completely different approach, see \cite{DSt31}.
\end{rk}

\end{document}